\documentclass[11pt]{amsart}
\allowdisplaybreaks

\usepackage{amsmath,amssymb,verbatim,graphicx,amsthm,color,url, enumerate}
\usepackage[top=0.5in, bottom=0.5in, left=1.2in, right=1.2in]{geometry}
\usepackage{cite}
\usepackage{color}
\usepackage{fancyhdr}
\usepackage{datetime}
\usepackage{csquotes}

\usepackage[backref=page,colorlinks=true]{hyperref} 
\usepackage{ifthen}
\usepackage{tikz}
\usetikzlibrary{decorations.pathreplacing}
\usepackage{tikz-cd}
\renewcommand*{\backref}[1]{}
\renewcommand*{\backrefalt}[4]{\tiny
	\ifcase #1 (\textbf{NOT CITED.})%
	\or    (Cited on page~#2.)%
	\else   (Cited on pages~#2.)%
	\fi}

\allowdisplaybreaks

\usepackage{blindtext}
\usepackage{scrextend}
\addtokomafont{labelinglabel}{\sffamily}

\newcommand{\vertiii}[1]{{\left\vert\kern-0.25ex\left\vert\kern-0.25ex\left\vert #1
    \right\vert\kern-0.25ex\right\vert\kern-0.25ex\right\vert}}
\newcommand{\braciii}[1]{{\left[\kern-0.25ex\left[\kern-0.25ex\left[ #1
		\right]\kern-0.25ex\right]\kern-0.25ex\right]}}

\newcommand{\ZZ}{\mathbb{Z}}

\newcommand{\NN}{\mathbb{N}}

\theoremstyle{plain}

\newtheorem{theorem}{Theorem}
\newtheorem*{theorem*}{Theorem}

\newtheorem{statement}[theorem]{Statement}

\theoremstyle{definition}

\title{A remark on "A non-singular dynamical system without maximal ergodic inequality" by E. H. El Abdalaoui}

\usepackage{setspace}

\author{Idris Assani}
\address{Department of Mathematics, The University of North Carolina at Chapel Hill}
\email{assani@email.unc.edu}
\urladdr{https://idrisassani.web.unc.edu/}

\begin{document}
	\maketitle
 \begin{abstract}
  In this note we would like to correct a comment made by E.H. El Abdaloui about my work.
 \end{abstract}
	In this note, I would like to address a comment regarding my work \cite{A1} that was made by El Houcein El Abdalaoui in his work, "A non-singular dynamical system without maximal ergodic inequality," which was published in this journal \cite{Ab19}. In particular, he writes that the strategy in my work "fails" without providing any evidence as to why that is the case.
	
	To provide some context, we summarize what was done in my paper. The goal of the preprint \cite{A1} was to show the pointwise convergence of multiple ergodic averages with several commuting transformations. More precisely:
	
	\begin{statement}\label{mea}
	 Let $k \in \NN$, let $(X, \mathcal{F}, \mu)$ be a probability measure space, let $T_1, T_2, \ldots, T_k: X \to X$ be $\mu$-invariant transformations on $X$ that are pairwise commuting (i.e. for every $i, j \in \{1, 2, \ldots, k\}$, $T_iT_j = T_jT_i$), and let $f_1, f_2, \ldots, f_k \in L^\infty(\mu)$. Then for $\mu$-a.e. $x \in X$, the limit
	\[\lim_{N \to \infty} \frac{1}{N} \sum_{n=1}^N \prod_{j=1}^k f_j(T^{jn}x) \]
	exists.
	\end{statement}
	Let $k \in \NN$, and $\Delta := \{z \in X^k: z = (x, x, \ldots, x)\}$ denote the diagonal of $X^k$. Furthermore, we denote $\mu_\Delta$ to be the diagonal measure of $X^k$, i.e. for any functions $g_1, g_2, \ldots, g_k \in L^\infty(\mu)$, we have
	$$
	\int g_1(x_1) g_2(x_2) \cdots g_k(x_k) \, d\mu_\Delta(x_1, x_2, \ldots,
	x_k) := \int g_1(x)g_2(x) \cdots g_k(x) \, d\mu(x).
	$$
	Statement \ref{mea} is equivalent of showing the existence of the limit
	\[ \lim_{N \to \infty} \frac{1}{N} \sum_{n=1}^N \prod_{j=1}^k f_j (T_j^{n}x_j)\]
	for $\mu_\Delta$-a.e. $(x_1, \ldots, x_k) \in X^k$.
	
	The issue with the diagonal measure is that measure is supported on $\Delta$, and does not see the action generated by the transformations $T_1 \times T_2 \times \cdots T_k$. Because of this, I introduced a
	non-singular dynamical system $(X^k, \mathcal{B}^k, \nu, \phi)$, which we refer to as "diagonal-orbit system" in my later work \cite{A2}, as follows: 
	We define $\phi := T_1 \times T_2
	\times \cdots \times T_k$, and $\nu$ a probability measure such that for any
	$A \in \mathcal{B}^k$,
	$$
	\nu(A) := \frac{1}{3} \sum_{n \in \ZZ} \frac{1}{2^{|n|}}\mu_\Delta
	(\phi^n A).
	$$
	
	While $\phi$ is not $\nu$-invariant, it is non-singular with respect to $\nu$. In fact, one of the crucial properties of the measure $\nu$ is that, for every $A \in \mathcal{B}^k$, we have
	\[ \nu(\phi^{-1}A) \leq 2\nu(A) \, . \]
	Furthermore, if $A \subset X^k$ for which $\nu(A) = 0$, then $\mu_\Delta(A) = 0$.
	
	In order to show that Statement \ref{mea} is true, it would be sufficient to prove the following: 
	
	\begin{statement}\label{mea-do}
	Let $k \in \NN$. $(X^k, \mathcal{B}^k, \nu, \phi := T_1 \times T_2 \times \cdots \times T_k)$ be the diagonal-orbit system of the system $(X, \mathcal{B}, \mu, T_1, T_2, \ldots, T_k)$, where $T_1, \ldots, T_k: X \to X$ are pairwise commuting measure-preserving transformations. Let $f_1, f_2, \ldots, f_k \in L^\infty(\mu)$, and let $F \in L^\infty(\nu)$ such that $F(x_1, x_2, \ldots, x_k) := \prod_{j=1}^k f_j(x_j)$. Then for $\nu$-a.e. $z \in X^k$, the limit
	\[ \lim_{N \to \infty} \frac{1}{N} \sum_{n=1}^N F(\phi^n z) \]
	exists.
	\end{statement}
	
	The following are quoted from p. 961 of the manuscript, paragraphs 2 (starting from line 5 from above), 3 (line 18 from above), and 6 (line 7 from below)  (the parts that have been omitted are denoted as ". . .", and the italicized words in a parenthesis (\textit{like this}) are my comments; the citation keys in the quote refers to the bibliography in Abdalaoui's paper): \\
	
	\begin{displayquote}
		"We recall that the $L^2$ convergence of the Furstenberg ergodic averages has been intensively studied and the topics is nowadays very rich. These studies were originated in the seminal work of Furstenberg on
		Szemeredi's theorem in [21]. Furstenberg--Katznelson--Ornstein in [22] proved that the $L^2$-norm convergence
		holds for $T_i = T^i$ and $T$ weakly mixing. Later, . . . (\textit{the remainder of this paragraph is devoted to a short survey of the study of $L^2$-convergence of the Furstenberg averages}) \\
		
		For the pointwise convergence, partial results were obtained in [18], [5], [6] and [30] under some ergodic
		and spectral assumptions. I. Assani in [7] purposes (\textit{sic; it is highly likely that he means "proposes"}) to solve the problem by noticing that the action of the maps $T_i$ induced a dynamical system $(X^k, \bigotimes_{i=1}^k \mathcal{B}, \nu, \phi)$ where . . . (\textit{explanation of the diagonal-orbit system, as well as some definitions such as a non-singular dynamical system, pushforward measure, etc. for the next few paragraphs}) \\
		
		Unfortunately as we shall see in section 5, the strategy of [7] (\textit{this is my preprint}) fails since the maximal ergodic inequality does not hold for the non-singular dynamical system $(X^k, \mathcal{B}^k, \nu, \phi)$."\\
	\end{displayquote}
	
	To explain why my strategy "fails", Abdalaoui shows that the pointwise ergodic theorem does not hold on the diagonal-orbit system, by showing the following in the proof of Theorem 2.1: Let $(X, \mathcal{B}, \mu, T_1, \ldots, T_k)$ be a measure-preserving system, where each $T_i$ is weakly mixing and pairwise commuting. Let $A \subset X^k$ such that $A := \bigcup_{n \in \ZZ}\phi^n(\Delta)$. If the pointwise ergodic theorem holds on $(X^k, \mathcal{B}^k, \nu, \phi)$, then for $\nu$-a.e. $z \in X^k$ we should have
	\[\lim_{N \to \infty} \frac{1}{N} \sum_{n=1}^N 1_A(\phi^n z) = \otimes_{j=1}^k \mu(A) = 0.  \]
	However, for any $N \in \NN$, we have
	\[ \int \frac{1}{N} \sum_{n=1}^N 1_A(\phi^n z) \, d\nu(z) = 1 \,  \]
	which is a contradiction. Hence, the pointwise ergodic theorem on the diagonal-orbit system does not hold.
	
	There are a few problems regarding this statement. First of all,\textit{ I have never been studying the pointwise ergodic theorem on the diagonal-orbit system}, meaning that I have never tried to prove the pointwise convergence of the averages \[\frac{1}{N} \sum_{n=1}^N F(\phi^n z)\] for functions in $L^p(\nu)$ for some $ 1\leq p\leq \infty.$ I have been using the diagonal-orbit system to study multiple ergodic averages, but this is a completely different problem from what Abdalaoui is discussing. In particular, Statement \ref{mea-do} states that the function $F$ in $L^\infty(\nu)$ must be written as a product of $L^\infty(\mu)$ functions from the base space $(X, \mathcal{B}, \mu)$. The example that Abdalaoui uses above, $1_A$, clearly does not belong in this class of functions.
	
	Secondly, the density argument that Abdalaoui uses is quite unclear. Let $\mathcal{L} \in L^\infty(\nu)$ be the finite linear combination of the set
	\[ \{f_1 \otimes f_2 \otimes \cdots f_k: f_1, f_2, \ldots, f_k \in L^\infty(\mu)\} \]
	By the Stone-Weierstrauss theorem, it is true that $\mathcal{L}$ is dense in the set of continuous functions $C(X^k)$, with respect to the supremum norm. However, $C(X^k)$ is \textit{not} dense in $L^\infty(\nu)$, with respect to the $L^\infty$-norm. Hence, it is unclear which density argument one could use to support Abdalaoui's claim. It should be noted that the existence of the maximal ergodic inequality would grant the approximation argument, but that is not the case for this system (as the whole point of this paper seems to be that such maximal ergodic inequality does not exist). Furthermore, the maximal ergodic inequality is not necessary to guarantee the pointwise convergence of the multiple ergodic averages in this problem, since the diagonal-orbit system $(X^k, \mathcal{B}^k, \nu, \phi)$ is not measure-preserving, and we are only interested in a subset of functions in $L^\infty(\nu)$. And in neither versions of my preprint, I have attempted to use the maximal ergodic theorem for the diagonal-orbit system (the first attempt has the maximal ergodic inequality for measure-preserving systems, which is classical, to connect the maximal function of a function and  its integral (see Lemma 2 in \cite{A1}) while the second attempt does not use any maximal ergodic inequality).
	
	Finally, I have been puzzled by the claim that a certain mathematical method "fails" to work. It is possible to rigorously prove that a mathematical statement is false by providing a counterexample. However, I am not aware of a way of demonstrating that a certain method would not work to prove a desired result. For instance, I have a feeling that the Galois theory may not be an applicable tool to show the pointwise convergence of multiple ergodic averages---but I do not know how to write a mathematical statement of such feeling and rigorously prove it (and it is possible that my feeling is wrong!).
	
	As a side note, I am currently working on revising the aforementioned preprint \cite{A1}, and hoping to announce it on arXiv sometimes in a near future. I was also able to prove a problem raised by J.-P. Conze regarding coboundary of multiple ergodic averages using the diagonal-orbit system, which can be seen in \cite{A2}.

\end{document}